\documentclass[11pt]{amsart}
\usepackage{amsmath,amssymb}
\newtheorem{theorem}{Theorem}[section]
\newtheorem{corollary}[theorem]{Corollary}

\theoremstyle{definition}
\newtheorem{definition}[theorem]{Definition}

\begin{document}

\title{Two-point functions and their applications in geometry}
\author{Simon Brendle}
\address{Department of Mathematics \\ Stanford University \\ Stanford, CA 94305}
\begin{abstract}
The maximum principle is one of the most important tools in the analysis of geometric partial differential equations. Traditionally, the maximum principle is applied to a scalar function defined on a manifold, but in recent years more sophisticated versions have emerged. One particularly interesting direction involves applying the maximum principle to functions that depend on a pair of points. This technique is particularly effective in the study of problems involving embedded surfaces. 

In this survey, we first describe some foundational results on curve shortening flow and mean curvature flow. We then describe Huisken's work on the curve shortening flow where the method of two-point functions was introduced. Finally, we discuss several recent applications of that technique. These include sharp estimates for mean curvature flow, as well as the proof of Lawson's 1970 conjecture concerning minimal tori in $S^3$.
\end{abstract}
\maketitle 

\section{Background on minimal surfaces and mean curvature flow}

Minimal surfaces are among the most important objects studied in differential geometry. A minimal surface is characterized by the fact that it is a critical point for the area functional; in other words, if we deform the surface while keeping the boundary fixed, then the surface area is unchanged to first order. Surfaces with this property serve as mathematical models for soap films in physics.

The first variation of the surface area can be expressed in terms of the curvature of the surface. To explain this, let us consider the most basic case of a surface $M$ in $\mathbb{R}^3$. The curvature of $M$ at a point $p \in M$ is described by a quadratic form defined on the tangent plane $T_p M$, which is referred to as the second fundamental form of $M$. The eigenvalues of the second fundamental form are the principal curvatures of $M$. We can think of the principal curvatures as follows: Given any point $p \in M$, we may locally represent the surface $M$ as a graph over the tangent plane $T_p M$. The principal curvatures of $M$ at $p$ can then be interpreted as the eigenvalues of the Hessian of the height function at the point $p$. For example, for a sphere of radius $r$, the principal curvatures are both equal to $\frac{1}{r}$; similarly, the principal curvatures of a cylinder of radius $r$ are equal to $\frac{1}{r}$ and $0$. If the principal curvatures of $M$ at the point $p$ have the same sign, then the surface $M$ will be convex or concave near $p$, depending on the orientation. On the other hand, if the principal curvatures at a point $p$ have opposite signs, then the surface $M$ will look like a saddle locally near $p$. The sum of the principal curvatures is referred to as the mean curvature of $M$ and is denoted by $H$.

It turns out that the mean curvature of a surface can be interpreted as the $L^2$-gradient of the area functional. This leads to the following equivalent characterization of minimal surfaces:

\begin{theorem} 
\label{characterization.of.minimal.surfaces}
Let $M$ be a two-dimensional surface in $\mathbb{R}^3$. Then the following statements are equivalent: 
\begin{itemize}
\item[(i)] $M$ is a minimal surface; that is, $M$ is a critical point of the area functional. 
\item[(ii)] The mean curvature of $M$ vanishes identically.
\item[(iii)] The restrictions of the coordinate functions in $\mathbb{R}^3$ are harmonic functions on $M$. In other words, if $F$ is a conformal parametrization of $M$, then the component functions $F_1,F_2,F_3$ are harmonic functions in the usual sense.
\end{itemize}
\end{theorem}

Theorem \ref{characterization.of.minimal.surfaces} is a classical fact in the theory of minimal surfaces. The standard examples of minimal surfaces in $\mathbb{R}^3$ include the catenoid 
\[\{(\cosh s \, \cos t,\cosh s \, \sin t,s): s,t \in \mathbb{R}\}\] 
and the helicoid 
\[\{(\sinh s \, \cos t,\sinh s \, \sin t,t): s,t \in \mathbb{R}\}.\] 
The global theory of minimal surfaces in $\mathbb{R}^3$ has attracted considerable interest in recent years. We will not discuss this theory here; instead, we refer the reader to the excellent texts \cite{Colding-Minicozzi1}, \cite{Colding-Minicozzi2}, \cite{Lawson-book}, \cite{Osserman}, and \cite{Simon2}.

We note that the definition of the principal curvatures and the mean curvature can be adapted to higher dimensional surfaces, as well as to surfaces in curved background manifolds. In particular, the notion of a minimal surface makes sense in this more general setting. This leads to some interesting new phenomena. For example, while there are no closed minimal surfaces in Euclidean space, it turns out that there are interesting examples of closed minimal surfaces in the unit sphere. The most basic examples of minimal surfaces in $S^3$ are the equator 
\[\{(x_1,x_2,x_3,x_4) \in \mathbb{R}^4: x_1^2+x_2^2+x_3^2=1, \, x_4=0\} \subset S^3.\] 
and the Clifford torus 
\[\Big \{ (x_1,x_2,x_3,x_4) \in \mathbb{R}^4: x_1^2+x_2^2=x_3^2+x_4^2=\frac{1}{2} \Big \} \subset S^3.\] 
Note that the principal curvatures of the equator both vanish. By contrast, the principal curvatures of the Clifford torus are equal to $1$ and $-1$.

One of the most important tools in the study of minimal surfaces is the maximum principle. Traditionally, one applies the maximum principle to some quantity involving the principal curvatures. One of the earliest examples of such an argument is the following theorem due to Simons: 

\begin{theorem}[J.~Simons \cite{Simons}]
\label{simons.thm}
Let $F: M \to S^{n+1}$ be an immersed minimal hypersurface such that $|A|^2 = \lambda_1^2+\hdots+\lambda_n^2 < n$, where $\lambda_1,\hdots,\lambda_n$ denote the principal curvatures. Then $F$ is a totally geodesic $n$-sphere.
\end{theorem}

Theorem \ref{simons.thm} can be viewed as a rigidity theorem for the equator in $S^{n+1}$. Note that the condition $|A|^2 < n$ is optimal; for example, the Clifford torus in $S^3$ satisfies $|A|^2=2$.

The proof of Simons' Theorem involves an application of the maximum principle to the function $|A|^2$. The Codazzi equations, together with the minimal surface equation, imply that 
\[\Delta h_{ik} + (|A|^2 - n) \, h_{ik} = 0,\] 
where $h_{ik}$ denotes the second fundamental form. From this, Simons deduced that 
\[\Delta (|A|^2) - 2 \, |\nabla A|^2 + 2 \, (|A|^2 - n) \, |A|^2 = 0.\] 
Consider now a point $p$ where the function $|A|^2$ attains its maximum. Clearly, $\Delta (|A|^2) \leq 0$ at $p$. This implies $(|A|^2 - n) \, |A|^2 \geq 0$ at $p$. Since $|A|^2 < n$ at $p$, we conclude that $|A|^2 = 0$ at $p$. Since the function $|A|^2$ attains its maximum at the point $p$, it follows that $|A|^2 = 0$ everywhere. In other words, the hypersurface is a totally geodesic $n$-sphere. \\

We note that minimal surfaces in the sphere are in one-to-one correspondence to minimal cones in Euclidean space. Indeed, using similar methods, Simons \cite{Simons} was able to prove a nonexistence theorem for stable minimal cones of dimension less than $7$. This result plays a central role in the regularity theory of minimal hypersurfaces (see e.g. \cite{Bombieri-DeGiorgi-Giusti}, \cite{Simon1}).

We now turn to parabolic PDEs. The study of nonlinear heat equations in differential geometry has a long history, going back to the work of Eells and Sampson \cite{Eells-Sampson} on the harmonic map heat flow and the fundamental work of Hamilton on the Ricci flow (cf. \cite{Hamilton1}, \cite{Hamilton2}). The most basic example of a geometric evolution equation is the curve shortening flow. This evolution equation was first studied by Gage and Hamilton \cite{Gage-Hamilton}, and can be viewed as a nonlinear heat flow for curves in the plane. We briefly review its definition. 

\begin{definition}
A family of immersed curves $F: S^1 \times [0,T) \to \mathbb{R}^2$ evolves by the curve shortening flow if 
\[\frac{\partial}{\partial t} F(x,t) = -\kappa(x,t) \, \nu(x,t),\] 
where $\nu$ denotes the unit normal to the curve $F(\cdot,t)$ and $\kappa$ denotes its geodesic curvature.
\end{definition}

\begin{theorem}[M.~Gage, R.~Hamilton \cite{Gage-Hamilton}]
\label{convergence.for.convex.curves}
Let $F: S^1 \times [0,T) \to \mathbb{R}^2$ be a convex solution of the curve shortening flow which is defined on a maximal time interval. Then, as $t$ approaches the singular time $T$, the curves $F(\cdot,t)$ shrink to a point, and converge to a circle after suitable rescaling.
\end{theorem}

There is a natural analogue of the curve shortening flow in higher dimensions, which is known as the mean curvature flow. This flow was first studied by Brakke \cite{Brakke} using techniques from geometric measure theory, and by Huisken \cite{Huisken1} using PDE techniques. 

\begin{definition}
Consider a manifold $M$ and a family of immersions $F: M \times [0,T) \to \mathbb{R}^2$. This family evolves by the mean shortening flow if 
\[\frac{\partial}{\partial t} F(x,t) = -H(x,t) \, \nu(x,t).\] 
Here, $\nu$ denotes the unit normal to the hypersurface $F(\cdot,t)$ and $H$ denotes the mean curvature (i.e. the sum of the principle curvatures).
\end{definition}

In the same way that minimal surface represent critical points of the area functional, the mean curvature flow has a natural geometric interpretation as the flow of steepest descent for the area functional. In the special case $n=1$, the mean curvature agrees with the geodesic curvature, and the mean curvature flow reduces to the curve shortening flow. The following result can be viewed as a higher-dimensional version of the theorem of Gage and Hamilton discussed above:

\begin{theorem}[G.~Huisken \cite{Huisken1}]
\label{mcf.convex}
Let $F: M \times [0,T) \to \mathbb{R}^{n+1}$, $n \geq 2$, be a convex solution of the mean curvature flow which is defined on a maximal time interval. Then, as $t$ approaches the singular time $T$, the hypersurfaces $F(\cdot,t)$ shrink to a point, and converge to a round sphere after suitable rescaling.
\end{theorem}

The starting point of Huisken's proof of Theorem \ref{mcf.convex} is a parabolic analogue of the Simons identity described above. More specifically, Huisken showed that the second fundamental form satisfies the evolution equation 
\[\frac{\partial}{\partial t} h_i^k = \Delta h_i^k + |A|^2 \, h_i^k.\] 
Hamilton's maximum principle for tensors (cf. \cite{Hamilton1}) then implies that the principal curvatures are bounded from above and below by fixed multiples of the mean curvature. In other words, one has $c \, H \leq \lambda_i \leq C \, H$, where $c$ and $C$ are positive constants that depend only on the initial hypersurface. Moreover, it turns out that the function 
\[f_\sigma = H^{\sigma-2} \, \Big ( |A|^2 - \frac{1}{n} \, H^2 \Big ).\] 
satisfies a differential inequality of the form 
\begin{align*} 
&\frac{\partial}{\partial t} f_\sigma - \Delta f_\sigma - 2 \, (1-\sigma) \, \Big \langle \frac{\nabla H}{H},\nabla f_\sigma \Big \rangle \\ 
&+ c^2 \, H^{\sigma-2} \, |\nabla H|^2 + \sigma \, |A|^2 \, f_\sigma \leq 0 
\end{align*} 
(cf. \cite{Huisken1}, Corollary 5.3). Huisken then used integral estimates and Stampacchia iteration to show that $f_\sigma \leq k$ for suitable constants $\sigma$ and $k$ that depends only on the initial hypersurface. This is a very sophisticated argument, which also makes use of the crucial Sobolev inequality due to Michael and Simon \cite{Michael-Simon}. This leads to the estimate 
\[|A|^2 - \frac{1}{n} \, H^2 \leq k \, H^{2-\sigma}.\] 
This inequality implies that any blow-up limit must satisfy $|A|^2 - \frac{1}{n} \, H^2 = 0$. Since $n \geq 2$, it follows that every blow-up limit is a round sphere by Schur's lemma.

\section{Curve shortening flow for embedded curves in the plane}

\label{curve.shortening}

In this section, we will consider solutions to the curve shortening flow that are embedded, i.e. free of self-intersections. Using the strict maximum principle, it is easy to see that embeddedness is preserved under the curve shortening flow. In other words, if $F_0$ is an embedding, then $F_t$ is an embedding for each $t \geq 0$. In 1997, Huisken obtained the following monotonicity formula for the curve shortening flow. This estimate is a quantitative version of the fact that embedded curves remain embedded under the evolution.

\begin{theorem}[G.~Huisken \cite{Huisken3}]
\label{chord.arc}
Let $F: S^1 \times [0,T) \to \mathbb{R}^2$ be a family of embedded curves in the plane which evolve by the curve shortening flow. Let $L(t)$ denote the curve of $F(\cdot,t)$, and let $d_t(x,y)$ denote the intrinsic distance of two points $x$ and $y$. Then the quantity
\[\sup_{x \neq y} \frac{L(t)}{|F(x,t) - F(y,t)|} \, \sin \frac{\pi \, d_t(x,y)}{L(t)}\] 
is monotone decreasing in $t$.
\end{theorem}

Let us sketch the proof of Theorem \ref{chord.arc}, following Huisken's paper \cite{Huisken3}. If the assertion is false, we can find times $t_0 < t_1$ and a real number $\alpha$ such that 
\begin{equation} 
\label{t_0}
\sup_{x \neq y} \frac{L(t_0)}{|F(x,t_0) - F(y,t_0)|} \, \sin \frac{\pi \, d_{t_0}(x,y)}{L(t_0)} < \alpha 
\end{equation}
and 
\begin{equation} 
\label{t_1}
\sup_{x \neq y} \frac{L(t_1)}{|F(x,t_1) - F(y,t_1)|} \, \sin \frac{\pi \, d_{t_1}(x,y)}{L(t_1)} > \alpha. 
\end{equation} 
Note that $\alpha > \pi$ in view of (\ref{t_0}). We now consider the function 
\[Z_\alpha(x,y,t) = \alpha \, |F(x,t) - F(y,t)| - L(t) \, \sin \frac{\pi \, d_t(x,y)}{L(t)}.\] 
In view of (\ref{t_0}) and (\ref{t_1}), there exists a time $\bar{t} \in (t_0,t_1)$ and a pair of points $\bar{x} \neq \bar{y}$ such that $Z_\alpha(\bar{x},\bar{y},\bar{t}) = 0$ and $Z_\alpha(x,y,t) \geq 0$ for all $x,y \in S^1$ and all $t \in (t_0,\bar{t})$.

Let $x$ be a local coordinate near $\bar{x}$, and let $y$ be a local coordinate near $\bar{y}$. We assume that $\frac{\partial F}{\partial x}(\bar{x},\bar{t})$ and $\frac{\partial F}{\partial y}(\bar{y},\bar{t})$ have unit length point. Moreover, we assume that $\frac{\partial}{\partial x}$ points away from $\bar{y}$, and $\frac{\partial}{\partial y}$ points away $\bar{x}$. Note that  
\[0 = \frac{\partial Z_\alpha}{\partial x}(\bar{x},\bar{y},\bar{t}) = \alpha \, \frac{\langle F(\bar{x},\bar{t})-F(\bar{y},\bar{t}),\frac{\partial F}{\partial x}(\bar{x},\bar{t}) \rangle}{|F(\bar{x},\bar{t})-F(\bar{y},\bar{t})|} - \pi \, \cos \frac{\pi \, d_{\bar{t}}(\bar{x},\bar{y})}{L(\bar{t})}\] 
and 
\[0 = \frac{\partial Z_\alpha}{\partial y}(\bar{x},\bar{y},\bar{t}) = -\alpha \, \frac{\langle F(\bar{x},\bar{t})-F(\bar{y},\bar{t}),\frac{\partial F}{\partial y}(\bar{y},\bar{t}) \rangle}{|F(\bar{x},\bar{t})-F(\bar{y},\bar{t})|} - \pi \, \cos \frac{\pi \, d_{\bar{t}}(\bar{x},\bar{y})}{L(\bar{t})}.\] 
In particular, subtracting the second identity from the first gives 
\[\Big \langle F(\bar{x},\bar{t})-F(\bar{y},\bar{t}),\frac{\partial F}{\partial x}(\bar{x},\bar{t}) \Big \rangle + \Big \langle F(\bar{x},\bar{t})-F(\bar{y},\bar{t}),\frac{\partial F}{\partial y}(\bar{y},\bar{t}) \Big \rangle = 0.\] 
This relation implies that the vector $\frac{\partial F}{\partial y}(\bar{y},\bar{t})$ is obtained from the vector $\frac{\partial F}{\partial y}(\bar{y},\bar{t})$ by reflection across the line orthogonal to $F(\bar{x},\bar{t})-F(\bar{y},\bar{t})$.

We now examine the second order partial derivatives of $Z_\alpha$. We compute
\begin{align*} 
&\frac{\partial^2 Z_\alpha}{\partial x^2}(\bar{x},\bar{y},\bar{t}) + \frac{\partial^2 Z_\alpha}{\partial y^2}(\bar{x},\bar{y},\bar{t}) + 2 \, \frac{\partial^2 Z_\alpha}{\partial x \, \partial y}(\bar{x},\bar{y},\bar{t}) \\ 
&= -\alpha \, \kappa(\bar{x},\bar{t}) \, \frac{\langle F(\bar{x},\bar{t})-F(\bar{y},\bar{t}),\nu(\bar{x},\bar{t}) \rangle}{|F(\bar{x},\bar{t})-F(\bar{y},\bar{t})|} \\ 
&+ \alpha \, \kappa(\bar{y},\bar{t}) \, \frac{\langle F(\bar{x},\bar{t})-F(\bar{y},\bar{t}),\nu(\bar{y},\bar{t}) \rangle}{|F(\bar{x},\bar{t})-F(\bar{y},\bar{t})|} \\ 
&+ \frac{4\pi^2}{L(\bar{t})} \, \sin \frac{\pi \, d_{\bar{t}}(\bar{x},\bar{y})}{L(\bar{t})}. 
\end{align*}
Finally, the time derivative of $Z_\alpha$ satisfies 
\begin{align*} 
\frac{\partial Z_\alpha}{\partial t}(\bar{x},\bar{y},\bar{t}) 
&= -\alpha \, \kappa(\bar{x},\bar{t}) \, \frac{\langle F(\bar{x},\bar{t})-F(\bar{y},\bar{t}),\nu(\bar{x},\bar{t}) \rangle}{|F(\bar{x},\bar{t})-F(\bar{y},\bar{t})|} \\ 
&+ \alpha \, \kappa(\bar{y},\bar{t}) \, \frac{\langle F(\bar{x},\bar{t})-F(\bar{y},\bar{t}),\nu(\bar{y},\bar{t}) \rangle}{|F(\bar{x},\bar{t})-F(\bar{y},\bar{t})|} \\ 
&+ \Big ( \sin \frac{\pi d_{\bar{t}}(\bar{x},\bar{y})}{L(\bar{t})} - \frac{\pi d_{\bar{t}}(\bar{x},\bar{y})}{L(\bar{t})} \, \cos \frac{\pi d_{\bar{t}}(\bar{x},\bar{y})}{L(\bar{t})} \Big ) \int_{S^1} \kappa^2 \\ 
&+ \pi \, \cos \frac{\pi d_{\bar{t}}(\bar{x},\bar{y})}{L(\bar{t})} \int_{\bar{x}}^{\bar{y}} \kappa^2. 
\end{align*} 
Since $\alpha > \pi$ and $Z_\alpha(\bar{x},\bar{y},\bar{t}) = 0$, the curve $F(\cdot,\bar{t})$ cannot have constant curvature. Using the Cauchy-Schwarz inequality and the Gauss-Bonnet theorem, we obtain 
\[\int_{S^1} \kappa^2 > \frac{1}{L(\bar{t})} \, \bigg ( \int_{S^1} \kappa \bigg )^2 = \frac{4\pi^2}{L(\bar{t})}.\] 
Similarly, 
\[\int_{\bar{x}}^{\bar{y}} \kappa^2 \geq \frac{1}{d_{\bar{t}}(\bar{x},\bar{y})} \, \bigg ( \int_{\bar{x}}^{\bar{y}} \kappa \bigg )^2 = \frac{4\xi^2}{d_{\bar{t}}(\bar{x},\bar{y})},\] 
where $\xi$ is defined by 
\[\cos \xi = \frac{\langle F(\bar{x},\bar{t})-F(\bar{y},\bar{t}),\frac{\partial F}{\partial x}(\bar{x},\bar{t}) \rangle}{|F(\bar{x},\bar{t})-F(\bar{y},\bar{t})|} = \frac{\pi}{\alpha} \, \cos \frac{\pi \, d_{\bar{t}}(\bar{x},\bar{y})}{L(\bar{t})}.\] 
Thus, 
\begin{align*} 
0 &\geq \frac{\partial Z_\alpha}{\partial t}(\bar{x},\bar{y},\bar{t}) - \frac{\partial^2 Z_\alpha}{\partial x^2}(\bar{x},\bar{y},\bar{t}) - \frac{\partial^2 Z_\alpha}{\partial y^2}(\bar{x},\bar{y},\bar{t}) - 2 \, \frac{\partial^2 Z_\alpha}{\partial x \, \partial y}(\bar{x},\bar{y},\bar{t}) \\ 
&> \frac{4\pi^2}{L(\bar{t})} \, \Big ( \sin \frac{\pi d_{\bar{t}}(\bar{x},\bar{y})}{L(\bar{t})} - \frac{\pi d_{\bar{t}}(\bar{x},\bar{y})}{L(\bar{t})} \, \cos \frac{\pi d_{\bar{t}}(\bar{x},\bar{y})}{L(\bar{t})} \Big ) \\ 
&+ \frac{4\xi^2}{d_{\bar{t}}(\bar{x},\bar{y})} \, \pi \, \cos \frac{\pi d_{\bar{t}}(\bar{x},\bar{y})}{L(\bar{t})} - \frac{4\pi^2}{L(\bar{t})} \, \sin \frac{\pi \, d_{\bar{t}}(\bar{x},\bar{y})}{L(\bar{t})} \\ 
&= \frac{4\pi}{d_{\bar{t}}(\bar{x},\bar{y})} \, \Big ( \xi^2 - \frac{\pi^2 d_{\bar{t}}(\bar{x},\bar{y})^2}{L(\bar{t})^2} \Big ) \, \cos \frac{\pi d_{\bar{t}}(\bar{x},\bar{y})}{L(\bar{t})}. 
\end{align*} 
On the other hand, the inequality $\alpha > \pi$ implies $\cos \xi \leq \cos \frac{\pi \, d_{\bar{t}}(\bar{x},\bar{y})}{L(\bar{t})}$, hence $\xi \geq \frac{\pi \, d_{\bar{t}}(\bar{x},\bar{y})}{L(\bar{t})}$. This is a contradiction. \\

Using Theorem \ref{chord.arc}, Huisken was able to show that the curve shortening flow shrinks every embedded curve in the plane to a round point. This was originally proved by Grayson \cite{Grayson} using different ideas (see also \cite{Hamilton3}):

\begin{theorem}[M.~Grayson \cite{Grayson}]
\label{convergence.for.embedded.curves}
Let $F: S^1 \times [0,T) \to \mathbb{R}^2$ be an embedded solution of the curve shortening flow which is defined on a maximal time interval. Then, as $t$ approaches the singular time $T$, the curves $F(\cdot,t)$ shrink to a point, and converge to a circle after suitable rescaling.
\end{theorem}

Theorem \ref{convergence.for.embedded.curves} can be deduced from Huisken's monotonicity formula above. Indeed, Theorem \ref{chord.arc} implies the chord-arc estimate 
\[\sin \frac{\pi \, d_t(x,y)}{L(t)} \leq C \, \frac{|F(x,t) - F(y,t)|}{L(t)}\] 
for some uniform constant $C$. In other words, the intrinsic distance between two points on the curve is bounded from above by a fixed multiple of the extrinsic distance. 

Now, if the flow does not converge to a circle after rescaling, then there exists a sequence of rescalings which converges to the grim reaper curve 
\[\Big \{ (\log \cos s,s): s \in (-\frac{\pi}{2},\frac{\pi}{2}) \Big \}\] 
(see \cite{Altschuler}, p.~492). An inspection of the grim reaper curve shows that we can find pairs of points with the property that their intrinsic distance is much larger than their extrinsic distance. This contradicts Huisken's chord-arc estimate. Thus, the grim reaper is ruled out as a blow-up limit, which implies that the flow converges to a circle after rescaling.

\section{Mean curvature flow for embedded mean convex hypersurfaces}

\label{mcf}

In this section, we consider solutions to the mean curvature flow in $\mathbb{R}^{n+1}$. As in the case of curves, the strict maximum principle implies that, if the initial hypersurface is embedded, then the flow will remain embedded for all subsequent times. Furthermore, if the initial hypersurface has positive mean curvature, then this will remain so for all time. 

In the following, we will consider a one-parameter family of embedded, mean-convex hypersurfaces $M_t$ in $\mathbb{R}^{n+1}$ which evolve by the mean curvature flow. In other words, we require that the hypersurfaces $M_t$ have positive mean curvature and are free of self-intersections. In a series of breathroughs \cite{White1}, \cite{White2}, \cite{White3}, \cite{White4}, White proved a number of deep results about the structure of the singular singularities that can arise along the flow. The following theorem is one of the central results of this theory: 

\begin{theorem}[B.~White \cite{White2}, p.~124] 
\label{mcf.singularities}
Let $M_t$, $t \in [0,T)$, be a family of smooth, embedded, mean convex hypersurfaces in $\mathbb{R}^{n+1}$ evolving under mean curvature flow. Given any sequence of numbers $t_i \in [0,T)$ and any sequence of points $p_i \in M_{t_i}$, there exists a subsequence which satisfies one of the following conditions: 
\begin{itemize} 
\item[(i)] There exists an open set $U \subset \mathbb{R}^{n+1}$ with the property that the hypersurfaces $M_{t_i} \cap U$ converge to a smooth, embedded hypersurface as $i \to \infty$. Moreover, the sequence $p_i$ converges to a point in $U$.
\item[(ii)] We have $H(p_i,t_i) \to \infty$, and the rescaled hypersurfaces $H(p_i,t_i) \, (M_{t_i} - p_i)$ converge in $C_{loc}^\infty$ to a smooth, embedded, and convex hypersurface. 
\end{itemize}
\end{theorem}

Theorem \ref{mcf.singularities} is in some ways analogous to Perelman's compactness theorem for ancient $\kappa$-solutions to the Ricci flow (cf. \cite{Perelman1}). The proof of Theorem \ref{mcf.singularities} is very subtle: a key step in the proof involves showing that static planes of multiplicity two or higher cannot arise as tangent flows (cf. \cite{White2}, Theorem 5). The proof of this fact relies on some fundamental new tools, such as the Sheeting Theorem and the Expanding Hole Theorem established in \cite{White1}. Another fundamental ingredient in the proof is Huisken's monotonicity formula for the mean curvature flow (cf. \cite{Huisken2}).

Theorem \ref{mcf.singularities} has several important consequences. For example, since the limiting surface in Theorem \ref{mcf.singularities} is smooth, White's theorem immediately implies a bound on the derivatives of the curvature. 

\begin{corollary}[B.~White \cite{White2}]
Let $M_t$, $t \in [0,T)$, be a family of smooth, embedded, mean convex hypersurfaces in $\mathbb{R}^{n+1}$ evolving under mean curvature flow. Then the derivatives of the curvature are bounded by $|\nabla A|^2 \leq C \, H^4$, where $C$ is a positive constant that depends only the initial hypersurface.
\end{corollary}

Another important consequence of Theorem \ref{mcf.singularities} is a lower bound for the inscribed radius. To describe this result, we first review the definition of the inscribed radius:

\begin{definition}
Given a closed hypersurface and a point $p$ on that hypersurface, the inscribed radius at $p$ is defined as the radius of the largest ball which is contained in the region enclosed by the hypersurface, and touches the hypersurface at $p$.
\end{definition}

The following result can be viewed as a quantitative version of the fact that the mean curvature flow preserves embeddedness:

\begin{theorem}[B.~White \cite{White2}; W.~Sheng, X.J.~Wang \cite{Sheng-Wang}; B.~Andrews \cite{Andrews}]
\label{inscribed.radius.non.sharp.bound}
Let $M_t$, $t \in [0,T)$, be a family of smooth, embedded, mean convex hypersurfaces in $\mathbb{R}^{n+1}$ evolving under mean curvature flow. Then the inscribed radius is bounded from below by $\frac{c}{H}$, where $c$ is a positive constant that depends only the initial hypersurface.
\end{theorem}

As mentioned above, Theorem \ref{inscribed.radius.non.sharp.bound} is a direct consequence of White's theorem on the structure of singularities (cf. Theorem \ref{mcf.singularities}). Indeed, if Theorem \ref{inscribed.radius.non.sharp.bound} is false, then we can find a sequence of 
numbers $t_i \in [0,T)$ and points $p_i \in M_{t_i}$ such that $H(p_i,t_i) \, \rho_i \to 0$, where $\rho_i$ denotes the 
inscribed radius of $M_{t_i}$ at the point $p_i$. We now apply Theorem \ref{mcf.singularities} to this sequence of points. If statement (i) in Theorem \ref{mcf.singularities} holds, then the inscribed radius of $M_{t_i}$ at the point $p_i$ is bounded away from zero, which contradicts the fact that $H(p_i,t_i) \, \rho_i \to 0$. On the other hand, if statement (ii) in Theorem \ref{mcf.singularities} holds, then the rescaled hypersurfaces $H(p_i,t_i) \, (M_{t_i} - p_i)$ converge in $C_{loc}^\infty$ to a smooth, convex, embedded hypersurface. In particular, the inscribed radius of $H(p_i,t_i) \, (M_{t_i} - p_i)$ at the origin must have a positive lower bound. This again contradicts the fact that $H(p_i,t_i) \, \rho_i \to 0$. 

Sheng and Wang \cite{Sheng-Wang} provided an alternative proof of Theorem \ref{inscribed.radius.non.sharp.bound}. Their proof again uses compactness results and contradiction arguments.

Andrews \cite{Andrews} recently gave another proof of Theorem \ref{inscribed.radius.non.sharp.bound} based on a direct monotonicity argument. To explain this argument, it is convenient to parametrize the surfaces $M_t$ by a map $F: M \times [0,T) \to \mathbb{R}^{n+1}$. Let $\mu(x,t)$ denote the reciprocal of the inscribed radius at the point $x$ at time $t$; that is, 
\[\mu(x,t) = \sup_{y \in M \setminus \{x\}} \frac{2 \, \langle F(x,t)-F(y,t),\nu(x,t) \rangle}{|F(x,t)-F(y,t)|^2}.\] 
Using a method similar in spirit to the one employed in Huisken's paper \cite{Huisken3}, Andrews was able to show that the function $\mu$ satisfies the inequality 
\begin{equation}
\label{zzz} 
\frac{\partial \mu}{\partial t} \leq \Delta \mu + |A|^2 \, \mu 
\end{equation} 
in the viscosity sense. On the other hand, it is well-known that the mean curvature satisfies the equation 
\[\frac{\partial H}{\partial t} = \Delta H + |A|^2 \, H.\] 
Thus, if $H$ is positive, the maximum principle implies that the ratio $\frac{\mu}{H}$ is uniformly bounded from above. In other words, the inscribed radius is bounded from below by $\frac{c}{H}$ for some small positive constant $c$. \\

Note that the estimate in Theorem \ref{inscribed.radius.non.sharp.bound} is not sharp. We next describe a sharp estimate for the inscribed radius of surfaces evolving by mean curvature flow. This estimate was established recently in \cite{Brendle5}. 

\begin{theorem}[S.~Brendle \cite{Brendle5}]
\label{inscribed.radius.sharp}
Let $F: M \times [0,T) \to \mathbb{R}^{n+1}$ be a family of embedded mean convex hypersurfaces which evolve by the mean curvature flow. Then, given any constant $\delta>0$, there exists a constant $C(\delta)$ with the property that $\mu \leq (1+\delta) \, H$ whenever $H \geq C(\delta)$. In other words, the inscribed radius is at least $\frac{1}{(1+\delta) \, H}$ at all points where the curvature is greater than $C(\delta)$.
\end{theorem}

Note that the estimate in Theorem \ref{inscribed.radius.sharp} is optimal. The estimate in Theorem \ref{inscribed.radius.sharp} can be generalized to solutions of the mean curvature flow in Riemannian manifolds (cf. \cite{Brendle6}), but we will focus here on the Euclidean case for simplicity. 

The proof of Theorem \ref{inscribed.radius.sharp} builds on Andrews' proof of Theorem \ref{inscribed.radius.non.sharp.bound} discussed above (cf. \cite{Andrews}), but several new ingredients are required in order to obtain a sharp estimate. First, we show that the function $\mu$ satisfies the differential inequality
\begin{equation} 
\label{ineq.a}
\frac{\partial \mu}{\partial t} \leq \Delta \mu + |A|^2 \, \mu - \sum_{i=1}^n \frac{2}{\mu-\lambda_i} \, (D_i \mu)^2,
\end{equation}
where the inequality is interpreted in the viscosity sense. Second, we show that, for each $t$, the function $\mu$ satisfies an inequality of the form 
\begin{align}
\label{ineq.b}
0 &\leq \Delta \mu + \frac{1}{2} \, |A|^2 \, \mu - \frac{1}{2} \, H \, \mu^2 + \frac{1}{2} \, n^3 \, (n\varepsilon \, \mu + K_1(\varepsilon)) \, \mu^2 \notag \\
&+ \sum_{i=1}^n \frac{1}{\mu-\lambda_i} \, D_i \mu \, D_i H \\ 
&+ \frac{1}{2} \, \big ( H + n^3 \, (n\varepsilon \, \mu + K_1(\varepsilon)) \big ) \, \sum_{i=1}^n \frac{1}{(\mu-\lambda_i)^2} \, (D_i \mu)^2 \notag
\end{align}
Here $\varepsilon$ is a positive real number which can be chosen arbitrarily small, and $K_1(\varepsilon)$ is a positive constant that depends on $\varepsilon$ and the initial hypersurface. The proof of the inequality (\ref{ineq.b}) does not directly use the fact that the surfaces evolve by mean curvature flow; it solely relies on the convexity estimates established by Huisken and Sinestrari \cite{Huisken-Sinestrari1}. 

Note that the inequality (\ref{ineq.a}) differs from the inequality (\ref{zzz}) by an extra gradient term which has a favorable sign. This term plays a crucial role in the argument. Roughly speaking, the inequality (\ref{ineq.a}) implies that the maximum of the function $\frac{\mu}{H}$ is strictly decreasing unless $\mu$ is constant. On the other hand, the relation (\ref{ineq.b}) tells us that $\mu$ cannot be constant unless $\frac{\mu}{H}$ is close to $1$ or smaller. In order to make this precise, we use integral estimates and Stampacchia iteration. This technique originated in Huisken's work on the mean curvature flow for convex hypersurfaces (see Theorem \ref{mcf.convex} above). Let us fix a positive number $\delta>0$, and let $\sigma \in (0,\frac{1}{2})$. The inequality (\ref{ineq.a}) implies that the function
\[f_\sigma = H^{\sigma-1} \, (\mu - (1+\delta) \, H)\]
satisfies the inequality 
\begin{align}
\label{evol.f}
&\frac{\partial}{\partial t} f_\sigma - \Delta f_\sigma - 2 \, (1-\sigma) \, \Big \langle \frac{\nabla H}{H},\nabla f_\sigma \Big \rangle \notag \\
&+ 2 \, \sum_{i=1}^n \frac{H^{\sigma-1}}{\mu-\lambda_i} \, (D_i \mu)^2 - \sigma \, |A|^2 \, f_\sigma \leq 0
\end{align}
Using (\ref{evol.f}) and the relation (\ref{ineq.b}), we can obtain an $L^p$-bound for the function $f_\sigma$. More precisely, we can find a positive constant $c_0$, depending only on $\delta$ and the initial hypersurface $M_0$, with the following property: if $p \geq \frac{1}{c_0}$ and $\sigma \leq c_0 \, p^{-\frac{1}{2}}$, then we have 
\[\int_{M_t} f_{\sigma,+}^p \leq C,\] 
where $C$ is a positive constant that depends only on $p$, $\sigma$, $\delta$, and the initial hypersurface $M_0$.

Having established $L^p$-bounds for the function $f_\sigma$, we can then use the Michael-Simon Sobolev inequality (cf. \cite{Michael-Simon}) and Stampacchia iteration to obtain a sup-bound for $f_\sigma$. More precisely, we can find positive numbers $\sigma$ and $k$, depending only on $\delta$ and the initial hypersurface, such that $f_\sigma \leq k$ everywhere. This means that 
\[\mu \leq (1+\delta) \, H + k \, H^{1-\sigma}.\] 
Since $\delta$ can be chosen arbitrarily small, this implies the desired estimate. \\

We note that Haslhofer and Kleiner \cite{Haslhofer-Kleiner2} have subsequently found an alternative proof of Theorem \ref{inscribed.radius.sharp}. Their argument uses a contradiction argument and relies on compactness theorems from \cite{Haslhofer-Kleiner1} (see also \cite{White1}, \cite{White2}).

In the remainder of this section, we discuss some applications of Theorem \ref{inscribed.radius.sharp}. The inscribed radius estimate in Theorem \ref{inscribed.radius.sharp} is most useful when $n=2$: in this case, Theorem \ref{inscribed.radius.sharp} precludes the formation of singularities of the form $\Gamma \times \mathbb{R}$, where $\Gamma$ is a curve with non-constant curvature. This can be viewed as an analogue of the cylindrical estimates of Huisken and Sinestrari (see \cite{Huisken-Sinestrari2}). 

In a joint work with Gerhard Huisken \cite{Brendle-Huisken}, we have recently defined a notion of mean curvature flow with surgery for mean convex surfaces in $\mathbb{R}^3$. 

\begin{theorem}[S.~Brendle, G.~Huisken \cite{Brendle-Huisken}]
\label{mcf.surgery}
Let $M_0$ be a closed, embedded surface in $\mathbb{R}^3$ with positive mean curvature. Then there exists a mean curvature flow with surgeries starting from $M_0$ which terminates after finitely many steps.
\end{theorem}

The construction in \cite{Brendle-Huisken} is based on the earlier work of Huisken and Sinestrari \cite{Huisken-Sinestrari2} in the higher-dimensional case, which in turn shares some common features with the surgery construction for the Ricci flow due to Hamilton \cite{Hamilton4}, \cite{Hamilton5} and Perelman \cite{Perelman1}, \cite{Perelman2}, \cite{Perelman3}. 

The main ingredients in the proof of Theorem \ref{mcf.surgery} are the convexity estimates of Huisken and Sinestrari \cite{Huisken-Sinestrari1}; the inscribed radius estimate in Theorem \ref{inscribed.radius.sharp}; a local gradient estimate due to Haslhofer and Kleiner \cite{Haslhofer-Kleiner1}; as well as a pseudolocality principle for mean curvature flow established in \cite{Brendle-Huisken}.

\section{Minimal surfaces in $S^3$ and Lawson's conjecture}

In this section, we discuss a conjecture due to Lawson concerning minimal surfaces in the three-dimensional unit sphere $S^3$. We first describe the background of this conjecture. In a landmark paper published in the late 1960s, Lawson \cite{Lawson2} constructed an infinite family of minimal tori in $S^3$ (see also \cite{Hsiang-Lawson}). These tori are immersed, but, with the exception of the Clifford torus, they fail to be embedded. In the same paper \cite{Lawson2}, Lawson also constructed an infinite family of embedded minimal surfaces of higher genus:

\begin{theorem}[H.B.~Lawson, Jr. \cite{Lawson2}]
There exists at least one embedded minimal surface in $S^3$ for any given genus. Moreover, there are at least two such surfaces unless the genus is a prime number.
\end{theorem}

The proof of Lawson's theorem uses the symmetries of $S^3$ in an extremely ingenious way. The key idea is to start from a least area disk spanning a geodesic quadrilateral. One then reflects this minimal disk across the edges of the geodesic quadrilateral. Repeating this process, one obtains a collection of minimal disks, each of which is contained in a geodesic tetrahedron, and these geodesic tetrahedra are mutually disjoint. The union of these minimal disks then gives a smooth embedded minimal surface in $S^3$ without boundary. 

In 1985, Karcher, Pinkall, and Sterling \cite{Karcher-Pinkall-Sterling} obtained additional examples of embedded minimal surfaces in $S^3$; their construction is closely related to Lawson's and involves the use of tesselations of $S^3$ into cells that have the symmetries of a Platonic solid in $\mathbb{R}^3$. Another related construction was found recently by Choe and Soret \cite{Choe-Soret}. Furthermore, Kapouleas and Yang \cite{Kapouleas-Yang} have constructed another infinite family of embedded minimal surfaces in $S^3$. This construction uses gluing techniques and the implicit function theorem; the resulting examples look like two nearby copies of the Clifford torus, which are joined by a large number of catenoid bridges. For related gluing constructions for minimal surfaces, we refer the reader to the survey paper \cite{Kapouleas-survey}, Section 2.4.

We now turn to uniqueness questions. In 1966, Almgren \cite{Almgren} proved that any minimal surface in $S^3$ of genus $0$ (i.e. which is homeomorphic to $S^2$) must be congruent to the equator in $S^3$. Almgren's proof uses the method of Hopf differentials. In view of the results above, it is natural to ask whether the Clifford torus is the only embedded minimal surface in $S^3$ of genus $1$. This was conjectured by Lawson \cite{Lawson3} in 1970, and confirmed by the author in \cite{Brendle1} (see also the survey paper \cite{Brendle4}).

\begin{theorem}[S.~Brendle \cite{Brendle1}]
\label{lawson.conj}
Let $F: \Sigma \to S^3$ be an embedded minimal torus in $S^3$. Then $F$ is congruent to the Clifford torus.
\end{theorem}

The proof of Theorem \ref{lawson.conj} relies on a sharp estimate for a two-point function which is obtained using the maximum principle. This estimate shares some common features with the results on curve shortening flow and mean curvature flow discussed in Sections \ref{curve.shortening} and \ref{mcf} (cf. \cite{Huisken3}, \cite{Andrews}, \cite{Brendle5}). However, there are a number of crucial differences. For example, we now have to deal with the curvature of the ambient space $S^3$. Furthermore, the arguments in \cite{Andrews} and \cite{Brendle5} both rely in a crucial way on the positivity of the mean curvature; in particular, these arguments do not work for minimal surfaces. 

Another central ingredient in the proof is a theorem due to Lawson \cite{Lawson2} which links the genus of a minimal surface in $S^3$ to the number of umbilical points. In particular, Lawson's result implies that a minimal surface in $S^3$ of genus $1$ has no umbilical points; in other words, the function $|A|$ is strictly positive everywhere.

We now sketch the proof of Theorem \ref{lawson.conj}. Suppose that $F: \Sigma \to S^3$ is an embedded minimal torus in $S^3$ which is not congruent to the Clifford torus, and let $\nu(x)$ denotes a unit normal vector field to the surface. In other words, $\nu(x)$ is tangential to $S^3$, but orthogonal to the tangent plane to the minimal surface. For each $\alpha \geq 1$, we consider the function 
\[Z_\alpha(x,y) = \frac{\alpha}{\sqrt{2}} \, |A(x)| \, (1 - \langle F(x),F(y) \rangle) + \langle \nu(x),F(y) \rangle.\] 
Since the surface is embedded and the function $|A|$ is positive everywhere, we conclude that $Z_\alpha(x,y) \geq 0$ for all $x,y \in \Sigma$, provided that $\alpha$ is sufficiently large. 

We claim that $Z_1(x,y) \geq 0$ for all $x,y \in \Sigma$. Suppose that this is false. Let us define 
\[\kappa = \inf \{\alpha \geq 1: \text{\rm $Z_\alpha(x,y) \geq 0$ for all $x,y \in \Sigma$}\} > 1.\] 
Clearly, $Z_\kappa(x,y) \geq 0$ for all points $x,y \in \Sigma$. Furthermore, we can find a pair of points $\bar{x} \neq \bar{y}$ such that $Z_\kappa(\bar{x},\bar{y}) = 0$. 

Since the function $Z_\kappa$ attains a local minimum at the point $(\bar{x},\bar{y})$, the first variation of $Z_\kappa$ at that point vanishes; that is, 
\begin{equation} 
\label{first.variation.x}
\frac{\partial Z_\kappa}{\partial x_i}(\bar{x},\bar{y}) = 0 
\end{equation}
and 
\begin{equation} 
\label{first.variation.y}
\frac{\partial Z_\kappa}{\partial y_i}(\bar{x},\bar{y}) = 0. 
\end{equation} 
The relation (\ref{first.variation.y}) provides an important piece of geometry information: it implies that the tangent plane $dF(T_{\bar{y}} \Sigma) \subset T_{F(\bar{y})} S^3$ is obtained from the tangent plane $dF(T_{\bar{x}} \Sigma) \subset T_{F(\bar{x})} S^3$ by reflection across the hyperplane orthogonal to $F(\bar{x})-F(\bar{y})$.

We next examine the second derivatives of $Z_\kappa$. Let $(x_1,x_2)$ and $(y_1,y_2)$ be geodesic normal coordinate systems around $\bar{x}$ and $\bar{y}$. One first computes the Laplacian of $Z_\kappa$ with respect to the $x$ variable. To that end, we use the Simons identity 
\[\Delta (|A|^2) - 2 \, |\nabla A|^2 + 2 \, (|A|^2 - 2) \, |A|^2 = 0.\] 
For a two-dimensional minimal surface, we have $|\nabla A|^2 = 2 \, \big | \nabla |A| \big |^2$. This implies 
\begin{equation} 
\label{simons.id}
\Delta |A| - \frac{\big | \nabla |A| \big |^2}{|A|} + (|A|^2 - 2) \, |A| = 0
\end{equation} 
(cf. \cite{Simons}). Unfortunately, the gradient term in (\ref{simons.id}) has an unfavorable sign. On the other hand, by exploiting the relations (\ref{first.variation.x}) and (\ref{first.variation.y}), we are able to extract a gradient term, which turns out to have a good sign. Furthermore, this good term is strong enough to absorb the bad term coming from the Simons identity! After a lengthy calculation, we arrive at the inequality 
\begin{align} 
\label{a}
\sum_{i=1}^2 \frac{\partial^2 Z_\kappa}{\partial x_i^2}(\bar{x},\bar{y}) 
&\leq \sqrt{2} \, \kappa \, |A(\bar{x})| \notag \\ 
&- \frac{\kappa^2-1}{\sqrt{2} \, \kappa} \, \frac{|A(\bar{x})|}{1-\langle F(\bar{x}),F(\bar{y}) \rangle} \sum_{i=1}^2 \Big \langle \frac{\partial F}{\partial x_i}(\bar{x}),F(\bar{y}) \Big \rangle^2. 
\end{align}
Moreover, we can compute the Laplacian of $Z_\kappa$ with respect to the $y$ variable. This gives 
\begin{equation} 
\label{b}
\sum_{i=1}^2 \frac{\partial^2 Z_\kappa}{\partial y_i^2}(\bar{x},\bar{y}) = \sqrt{2} \, \kappa \, |A(\bar{x})|. 
\end{equation}
Finally, one has to exploit the mixed partial derivatives of $Z_\kappa$. We can arrange that the tangent vector $\frac{\partial F}{\partial y_i}(\bar{y}) \in T_{F(\bar{y})} S^3$ is obtained from the tangent vector $\frac{\partial F}{\partial x_i}(\bar{x}) \in T_{F(\bar{x})} S^3$ by reflection across the hyperplane orthogonal to $F(\bar{x})-F(\bar{y})$. With this choice of the coordinate system, we obtain 
\begin{equation} 
\label{c}
\sum_{i=1}^2 \frac{\partial^2 Z_\kappa}{\partial x_i \, \partial y_i} = -\sqrt{2} \, \kappa \, |A(\bar{x})|. 
\end{equation}
The proof of (\ref{c}) also uses the identity (\ref{first.variation.x}) above.

Combining (\ref{a}), (\ref{b}), and (\ref{c}), we obtain the differential inequality 
\begin{align} 
\label{main.inequality}
&\sum_{i=1}^2 \frac{\partial^2 Z}{\partial x_i^2}(\bar{x},\bar{y}) + 2 \sum_{i=1}^2 \frac{\partial^2 Z}{\partial x_i \, \partial y_i}(\bar{x},\bar{y}) + \sum_{i=1}^2 \frac{\partial^2 Z}{\partial y_i^2}(\bar{x},\bar{y}) \notag \\ 
&\leq -\frac{\kappa^2-1}{\sqrt{2} \, \kappa} \, \frac{|A(\bar{x})|}{1-\langle F(\bar{x}),F(\bar{y}) \rangle} \sum_{i=1}^2 \Big \langle \frac{\partial F}{\partial x_i}(\bar{x}),F(\bar{y}) \Big \rangle^2 \leq 0. 
\end{align} 
Note that (\ref{main.inequality}) is not a strict inequality, so we have not yet arrived at a contradiction. In order to finish the argument, we observe that there is a variant of (\ref{main.inequality}) which holds for all pairs of points $\bar{x} \neq \bar{y}$, not just at the minimum. Indeed, by adapting the proof of (\ref{main.inequality}), we can show that 
\begin{align} 
\label{pde}
&\sum_{i=1}^2 \frac{\partial^2 Z}{\partial x_i^2}(\bar{x},\bar{y}) + 2 \sum_{i=1}^2 \frac{\partial^2 Z}{\partial x_i \, \partial y_i}(\bar{x},\bar{y}) + \sum_{i=1}^2 \frac{\partial^2 Z}{\partial y_i^2}(\bar{x},\bar{y}) \notag \\ 
&\leq -\frac{\kappa^2-1}{\sqrt{2} \, \kappa} \, \frac{|A(\bar{x})|}{1-\langle F(\bar{x}),F(\bar{y}) \rangle} \sum_{i=1}^2 \Big \langle \frac{\partial F}{\partial x_i}(\bar{x}),F(\bar{y}) \Big \rangle^2 \\ 
&+ \Lambda(\bar{x},\bar{y}) \, \bigg ( Z_\kappa(\bar{x},\bar{y}) + \sum_{i=1}^2 \Big | \frac{\partial Z_\kappa}{\partial x_i}(\bar{x},\bar{y}) \Big | + \sum_{i=1}^2 \Big | \frac{\partial Z_\kappa}{\partial y_i}(\bar{x},\bar{y}) \Big | \bigg ) \notag
\end{align} 
for all pairs of points $\bar{x} \neq \bar{y}$. Here, $\Lambda$ is a positive function on the set $\{(x,y) \in \Sigma \times \Sigma: x \neq y\}$ which is bounded away from the diagonal. In view of Bony's version of the strict maximum principle for degenerate elliptic equations (cf. \cite{Bony}), we can conclude that the set 
\[\{x \in \Sigma: \text{\rm $Z_\kappa(x,y)=0$ for some $y \in \Sigma \setminus \{x\}$}\}\] 
is open and non-empty. This implies that $|A|$ is constant on a non-empty open set. By a theorem of Lawson \cite{Lawson1}, this can only happen if the surface is congruent to the Clifford torus, which contradicts our assumption. 

Thus, we must have $Z_1(x,y) \geq 0$ for all $x,y \in \Sigma$. In other words, we have 
\begin{equation} 
\label{z1}
\frac{1}{\sqrt{2}} \, |A(x)| \, (1 - \langle F(x),F(y) \rangle) + \langle \nu(x),F(y) \rangle \geq 0 
\end{equation}
for all $x,y \in \Sigma$. An analogous argument with $\nu$ replaced by $-\nu$ imply that 
\begin{equation} 
\label{z1.modified}
\frac{1}{\sqrt{2}} \, |A(x)| \, (1 - \langle F(x),F(y) \rangle) - \langle \nu(x),F(y) \rangle \geq 0 
\end{equation}
for all $x,y \in \Sigma$. By performing a Taylor expansion of the left hand side in (\ref{z1}) and (\ref{z1.modified}) for $y$ close to $x$, one can show that $\nabla |A|=0$ for any given point $x$. This implies that our surface is congruent to the Clifford torus, contrary to our assumption. This completes our sketch of the proof of Theorem \ref{lawson.conj}. \\

We note that the proof of Theorem \ref{lawson.conj} can be extended in various ways. For example, the arguments in \cite{Brendle1} can be adapted to constant mean curvature surfaces (see \cite{Andrews-Li}), as well as to surfaces that are immersed in the sense of Alexandrov (cf. \cite{Brendle2}). This leads to the result that every constant mean curvature torus in $S^3$ which is immersed in the sense of Alexandrov must be rotationally symmetric; see \cite{Andrews-Li} and \cite{Brendle2} for details. Moreover, the classification of constant mean curvature tori in $S^3$ with rotational symmetry can be reduced to the analysis of an ODE (see \cite{Perdomo}). It turns out that the class of surfaces which are immersed in the sense of Alexandrov is quite natural in this context; in particular, there is a large class of examples which are immersed in the sense of Alexandrov, but fail to be embedded. Finally, the proof of Theorem \ref{lawson.conj} can be extended to a class of Weingarten tori in $S^3$ (see \cite{Brendle3}).

\end{document}